\documentclass [a4paper,12pt]{article}
\usepackage[all]{xy}
\usepackage[english]{babel}
\usepackage[latin1]{inputenc}
\usepackage{amsmath}
\usepackage{amsthm}
\usepackage{amsfonts}
\usepackage{amssymb}
\usepackage{mathrsfs}
\usepackage{graphics}
\usepackage{graphicx}
\usepackage[lite]{amsrefs}

\newtheorem{teorema}{Theorem}[section]
\newtheorem{lemma}[teorema]{Lemma}
\newtheorem{propos}[teorema]{Proposition}
\newtheorem{corol}[teorema]{Corollary}
\theoremstyle{remark}
\newtheorem{rem}{Remark}[section]
\newtheorem{ex}{Example}[section]
\theoremstyle{definition}
\newtheorem{defin}{Definition}[section]

\def\demo{\par\noindent{\bf Proof.\ }}
\def\enddemo{\ $\hfill\Box$\par\vskip.6truecm}

\def\de{\partial}
\def\debar{\overline{\de}}

\def\R{{\mathbb R}}
\def\N{{\mathbb N}}
\def\C{{\mathbb C}}

\def\Ol{{\mathcal O}}

\def\E{{\mathscr E}}
\def\Di{{\mathcal D}}
\def\Dii{{\mathscr D}}
\def\Ci{{\mathcal C}}
\def\Leb{{\mathcal L}}
\def\supp{\mathrm{supp}\,}
\def\loc{\mathrm{loc}}
\def\Lip{\mathrm{Lip}}
\def\rg{\mathrm{reg}}
\def\sg{\mathrm{sing}}
\def\H{{\mathcal H}}
\def\mass{\mathbf{M}}
\def\q{{\mathrm{q-loc}}}

\newcommand{\pair}[1]{\left\langle #1 \right\rangle}

\begin {document}
\title{Some applications of metric currents to complex analysis}
\author{Samuele Mongodi}
    
\maketitle

\begin{abstract}
The aim of this paper is to extend the theory of metric currents, developed by Ambrosio and Kirchheim, to complex spaces. We define the bidimension of a metric current on a complex space and we discuss the Cauchy-Riemann equation on a particular class of singular spaces. As another application, we investigate the Cauchy-Riemann equation on complex Banach spaces, by means of a homotopy formula.
\end{abstract}

%\tableofcontents

\section{Introduction}

In 2000, L. Ambrosio and B. Kirchheim introduced a new and powerful tool to extend some of the finer analytic techniques of geometric measure theory to the general setting of metric spaces; the theory of metric currents, which they develop in \cite{AK}, was later expanded by U. Lang, in \cite{Lang}, to include, at least in locally compact spaces, also the study of local phenomena.

Here we present some applications of this theory to complex analysis, tackling the problem of the Cauchy-Riemann equation on singular complex 
spaces and in complex Banach spaces. 

\medskip

The study of $\debar$ equation on singular spaces dates back to Fornaess and Gavosto (\cite{FornGav}) and Henkin and Polyakov (\cite{HenPol}), with different approaches and statements of the problem; in the last 20 years, the \emph{implicit} approach of Fornaess and Gavosto has been pursued by Fornaess, Vassiliadou, Ovrelid, Ruppenthal in a number of works (e.g. \cite{ForVas}, \cite{VasOv}, \cite{ForVas2}, \cite{Rupp} and \cite{Rupp1}, which we refer to for a more accurate bibliography), employing subtle $L^2$ techniques and an accurate analysis of the desingularization procedure for complex spaces. 

\medskip

It is worth mentioning that the representation formulae techniques used by Henking and Polyakov have been greatly generalized and improved in the recent work of Andersson and Samuelsson (\cite{And1}), leading to a general representation kernel for the $\debar$ equation on complex spaces.

\medskip

Complex analysis in infinite dimension, complex Banach  spaces and complex Banach manifolds are widely studied themes, starting from the work of Lelong on topological vector spaces, the study of the space of moduli of complex subspaces by Douady (\cite{Douady}), going on with the work of Dineen, Nachbin, Noverraz, Rickart and others (see \cite{Noverraz} or \cite{Dineen} for an extensive bibliography). However, the problem of the Cauchy-Riemann equation has been approached, to our best knowledge, only recently by Lempert (see \cites{Lemp1, Lemp2, Lemp3}).

\bigskip

The language of metric currents gives us the possibility to describe an intrinsic object, both on a complex singular space and in a complex Banach space, allowing us to exploit the peculiarities of its geometry to produce a solution to the Cauchy-Riemann equation.

The aim of this paper is to show these two applications of metric currents to complex analysis.

In Section 2 we recall the basic definitions both from \cite{AK} and \cite{Lang} and we give a detailed proof of the comparison theorem between metric currents and classical ones on a manifold. In Section 3 we introduce the concept of bidimension for a metric current on a finite dimensional space, showing that the usual properties of $(p,q)-$currents still hold, except for the existence of a Dolbeault decomposition.

Section 4 is devoted to the analysis of a particular class of complex spaces, whose geometry allows us to give a structure theorem for currents, solve the Cauchy-Riemann equation and characterize holomorphic currents.

In Section 5, we introduce the concept of bidimension of (global) metric currents on a Banach space and relate it to the behavior of the finite dimensional projections of the currents. In Section 6 we define a new class of currents, the \emph{quasi-local} metric currents, which are usual metric currents when restricted to bounded sets, and we give a definition of $(p,q)-$currents in this new class.
The last section shows how to employ these newly defined quasi-local currents in order to obtain a solution to the equation $\debar U=T$, when $T$ is of bidimension $(0,q)$ and its support is bounded; finally, we extend the result to a current $T$ with generic bidimension, $\debar-$closed, with bounded support.

\noindent{\bf{Acknowledgement}} I'm greatly indebted to my advisor, prof. G. Tomassini, for the advice he has always given me and the encouragement without which this article would never have been written. I would also like to thank prof. L. Ambrosio for his willingness to read the drafts and discuss them with me. 

I acknowledge the support of the ERC ADG GeMeThNES.

\section{Basic Notions}

Let $X$ be a complex space, possibly singular, with a distance $d$, compatible with its complex space topology. This distance can be the one associated to a K\"ahler structure on $X$, the Kobayashi distance (in case $X$ is Kobayashi-hyperbolic), or the distance induced by an embedding of $X$ into some hermitian complex manifold (e.g. $\C^n$).

\medskip

In case $X$ is locally compact, which is the case for finite dimensional complex spaces, we will denote by $\Di_k(X)$ the space of $k-$dimensional local metric currents, as defined in \cite{Lang}, which we refer to for the basic properties and definitions regarding such currents; we will also denote by $M_{k,\loc}(X)$ and $N_{k,\loc}(X)$ the spaces of local currents with locally finite mass and of locally normal currents, while $M_k(X)$ and $N_k(X)$ will be used for currents with finite mass and normal currents.

If $X$ is not locally compact, but it is complete, we will define only the spaces $M_k(X)$ and $N_k(X)$ following \cite{AK}.

\medskip

Following \cite{Lang}, we set
$$\Di(X)=\bigcup_{\substack{K\Subset X\\L>0}}\Lip_{K,L}=\bigcup_{\substack{K\Subset X\\L>0}}\{f\in\Lip_L(X)\ :\ \supp(f)\subset K\}$$
$$\Di^0(X)=\Di(X)\qquad \Di^m(X)=\Di(X)\times[\Lip_\loc(X)]^m$$
and we introduce the notation
$$\E^0(X)=\Lip_b(X)\qquad \E^m(X)=\Lip_b(X)\times[\Lip(X)]^m$$
where $\Lip_b$ is the algebra of bounded Lipschitz functions. We remark that this notations differs from the one used in \cite{AK}. The elements of this sets are called metric forms; by a slight abuse of notation, we will identify the $k+1-$tuple $(f,\pi_1,\ldots,\pi_k)$ with the expression $fd\pi_1\wedge\ldots\wedge d\pi_k$, in view of the product rule and the alternating property of metric currents.

\medskip

We recall the definition of a local $k-$dimensional metric current.

\begin{defin}A \emph{local $k-$dimensional metric current} is a functional $T:\Di^k(X)\to\C$ satisfying the following
\begin{enumerate}
\item $T$ is multilinear
\item whenever $(f^i,\pi^i)\to (f,\pi)$ pointwise, with locally uniformly bounded Lipschitz constants, $T(f^i,\pi^i)\to T(f,\pi)$
\item $T(f,\pi)=0$, whenever there is an index $j$ for which $\pi_j$ is constant on a neighborhood of $\supp f$.
\end{enumerate}
\end{defin}

We say that a local current $T$ is \emph{of finite mass} if there exists a finite Radon measure $\mu$ on $X$ such that
$$|T(f,\pi)|\leq \prod_{j=1}^k\Lip(\pi_j)\int_{X}|f|d\mu\qquad \forall (f,\pi)\in\Di^k(X)\;;$$
the infimum of such measures $\mu$ is called the \emph{mass} (measure) of $T$.
Requiring the estimate to hold only for the metric forms with support in a given compact $K\Subset X$ and allowing $\mu=\mu_K$ to depend on that compact, we get the condition of \emph{locally finite mass}.

The mass induces a norm on the currents with finite mass, by setting $\|T\|=\mu(X)$, where $\mu$ is the mass measure of $T$.

For sake of completeness, we recall also the definition of metric currents following  \cite{AK}.

\begin{defin} A \emph{$k-$dimensional metric current} is a functional $T:\E^k(X)\to \C$ satisfying the following
\begin{enumerate}
\item $T$ is multilinear
\item whenever $(f^i,\pi^i)\to (f,\pi)$ pointwise, with uniformly bounded Lipschitz constants, $T(f^i,\pi^i)\to T(f,\pi)$
\item $T(f,\pi)=0$, whenever there is an index $j$ for which $\pi_j$ is constant on a neighborhood of $\supp f$
\item $T$ has finite mass.
\end{enumerate}
\end{defin}

If the metric space isn't locally compact, we can still define local currents of finite mass with locally compact support; the space of $k-$dimensional currents in the sense of Ambrosio and Kirchheim is the completion of such a space of local $k-$dimensional currents of finite mass with respect to the mass norm (see the end of Section 4 in \cite{Lang}).

\medskip

Given a $k-$dimensional (local) metric current $T$, we will denote its boundary by $dT$ and, for any $\phi:X\to Y$ which is (locally) Lipschitz (and proper), $\phi_\sharp T$ will denote the pushforward of $T$ through $\phi$. If $T\in M_k(X)$ and $udv\in \E^h(X)$, $T\llcorner(udv)$ is the $(k-h)-$metric current obtained contracting $T$ with $udv$; if $T\in\Di_k(X$), we can define $T\llcorner (udv)$ for any $(u,v)\in\Lip_\loc(X)^{h+1}$. Sometimes, we will use the notation $T\llcorner (u,v)$ as an alternative to $T\llcorner (udv)$.

\medskip

Given a sequence of metric currents $T_j\in \Di_m(X)$, we say that $T_j\to T\in \Di_m(X)$ weakly if $T_j(f,\pi)\to T(f,\pi)$ for every $(f,\pi)\in\Di^m(X)$.

\subsection{Comparison theorems}

We denote by $\Dii$, $\mathfrak{N}$, $\mathfrak{F}$ the spaces of classical currents, normal classical currents and flat classical currents, as defined in \cite{Federer}; we denote by $\mass(T)$ the classical mass of a current.

For the proof of the following theorem, see \cite[Theorem 5.5]{Lang} and \cite[Theorem 11.1]{AK}.

\begin{teorema} \label{teo_comp}Let $U\subset\C^N$ be an open set, $N\geq1$. For every $m\geq0$ there exists an injective linear map $C_m:\Di_m(U)\to\Dii_m(U)$ such that
$$C_m(T)(fdg_1\wedge\ldots\wedge dg_m)=T(f,g_1,\ldots, g_m)$$
for all $(f,g_1,\ldots, g_m)\in\Ci^\infty_c(U)\times[\Ci^\infty(U)]^m$. The following properties hold:
\begin{enumerate}
\item for $m\geq1$, $d\circ C_m=C_{m-1}\circ d$;
\item for all $T\in\Di_m(U)$, $\|T\|\leq\mass(C_m(T))\leq{N\choose m}\|T\|$;
\item the restriction of $C_m$ to $N_{m,\loc}$ is an isomorphism onto $\mathfrak{N}_{m,\loc}$;
\item the image of $C_m$ contains the space $\mathfrak{F}_{m,\loc}(U)$.\end{enumerate}
\end{teorema}

We have an analogous result for manifolds.

\begin{teorema}\label{teo_comp_man}Let $U$ be an $N-$dimensional complex manifold, $N\geq 1$. For every $m\geq0$ there exists an injective linear map $C_m:\Di_m(U)\to \Dii_m(U)$ such that
$$C_m(T)(fdg_1\wedge\ldots\wedge dg_m)=T(f,g_1,\ldots, g_m)$$
for all $(f,g_1,\ldots, g_m)\in\Ci^\infty_c(U)\times[\Ci^\infty(U)]^m$. The following properties hold:
\begin{enumerate}
\item for $m\geq1$, $d\circ C_m=C_{m-1}\circ d$;
\item there exists a positive constant $c_1$ such that, for all $T\in\Di_m(U)$, $c_1^{-2}\|T\|\leq\mass(C_m(T))\leq c_1^{2}{N\choose m}\|T\|$;
\item the restriction of $C_m$ to $N_{m,\loc}$ is an isomorphism onto $\mathfrak{N}_{m,\loc}$;
\item the image of $C_m$ contains the space $\mathfrak{F}_{m,\loc}(U)$.\end{enumerate}
\end{teorema}
\demo There exists a locally finite cover $\{U_j\}_{j\in\N}$ of relatively compact open sets with bi-Lipschitz coordinate charts $\phi_j:U_j\to\Omega_j\subseteq\C^N$.

Let $\{\rho_j\}_{j\in\N}$ be a smooth partition of unity subordinated to the cover $\{U_j\}_{j\in\N}$; for each $j\in\N$, the current $T_j=T\llcorner\rho_j$ is supported in $U_j$, therefore belongs to $\Di_m(U_j)$.

The induced map $(\phi_j)_\sharp$ is an isomorphism between $\Di_m(U_j)$ and $\Di_m(\Omega_j)$; let $C_m^j$ be the linear injective map given by Theorem \ref{teo_comp} between $\Di_m(\Omega_j)$ and $\Dii_m(\Omega_j)$. The map $(\phi_j^{-1})_*$ is an isomorphism between $\Dii_m(\Omega_j)$ and $\Dii_m(U_j)$, which can be injected into $\Dii_m(U)$.

Therefore, for every $j\in\N$, we have the map 
$$T\mapsto T_j\mapsto (\phi_j)_\sharp T_j\mapsto C_m^j((\phi_j)_\sharp T_j)\mapsto (\phi_j^{-1})_*C_m^j((\phi_j)_\sharp T_j)=R_m^j(T)\;.$$
All the intermediate steps are linear, so the result is linear. We set
$$C_m(T)=\sum_j R_m^j(T)\;,$$
which is well defined because the cover is locally finite; moreover, the map $C_m$ is linear and injective, because, if $T\neq S$, then there exists $j$ such that $T_j\neq S_j$, so $R_m^j(T)\neq R_m^j(S)$ and we can assume that there exists an open subset of $U_j$ not contained in any other $U_k$.

Given $(f,g_1,\ldots, g_m)\in\Ci^\infty_c(U)\times[\Ci^\infty(U)]^m$, we have that $T(f,g)$ can be written as
$$\sum_j T(\rho_j f, g)=\sum_j (\phi_j)_\sharp T_j ((\rho_j\circ\phi_j^{-1})\cdot(f\circ\phi_j^{-1}), g_1\circ\phi_j^{-1},\ldots, g_m\circ\phi_j^{-1})$$
where the sums are indeed finite, because $f$ has compact support. Now, by Theorem \ref{teo_comp}, we have that the last sum is equal to
$$\sum_j C_m^j((\phi_j)_\sharp T_j)((\rho_j f)\circ\phi_j^{-1}dg_1\circ\phi_j^{-1}\wedge\ldots\wedge dg_m\circ\phi_j^{-1})=$$
$$\sum_j(\phi_j^{-1})_*C_m^j((\phi_j)_\sharp T_j)(\rho_j fdg_1\wedge\ldots \wedge dg_m)=\sum R_m^j(T)(fdg_1\wedge\ldots \wedge dg_m)$$
$$=C_m(fdg_1\wedge\ldots \wedge dg_m)\;.$$

\medskip

The conclusions of Theorem \ref{teo_comp} hold for $C_m^j$ and for the pushforward maps. We need to check that they still hold after contraction with a $0-$form and after a locally finite sum.

\smallskip

\noindent{\emph{(1)}} We know that, for $m\geq1$, $d\circ C_m^j=C_{m-1}^j\circ d$. We also have
$$d(T_j)=d(T_\llcorner\rho_j)=dT\llcorner\rho_j-T\llcorner(\sigma_j,\rho_j)=(dT)_j-T\llcorner(\sigma_j,\rho_j)\;,$$
with $\sigma_j$ a compactly supported smooth function equal to $1$ on $\supp\rho_j$ and to $0$ outside $U_j$, so
$$R^j_{m-1}(dT)=$$ 
$$(\phi_j^{-1})_*C_{m-1}^j((\phi_j)_\sharp (dT)_j)=(\phi_j^{-1})_*C_{m-1}^j((\phi_j)_\sharp (d(T_j)-T\llcorner(\sigma_j,\rho_j)))=$$
$$(\phi^{-1}_j)_*C_{m-1}^j(d(\phi_j)_\sharp T_j-((\phi_j)_\sharp T\llcorner(\sigma_j,\rho_j)))$$
$$=(\phi_j^{-1})_*(dC_m^j((\phi_j)_\sharp T_j)-C_{m-1}^j((\phi_j)_\sharp T\llcorner (\sigma_j,\rho_j)))=$$
$$dR_m^j(T)-S_{m-1}^j(T)\;.$$
For a given classical form $fdg_1\wedge\ldots\wedge dg_{m-1}$, we have that 
$$S_{m-1}^j(T)(fdg_1\wedge\ldots\wedge dg_{m-1})\neq0$$ 
only for a finite number of $j\in\N$ (namely, those such that $fd\rho_j \neq0$), so
$$\sum_{j:fd\rho_j\neq0}\!\!\!\!S_{m-1}^j(T)(fdg_1\wedge\ldots\wedge dg_{m-1})=$$
$$\sum(\phi_j^{-1})_* C_{m-1}^j((\phi_j)_\sharp T\llcorner (\sigma_j,\rho_j)))(fdg_1\wedge\ldots\wedge dg_{m-1})=$$
$$\sum T\llcorner (\sigma_j,\rho_j)(f,g_1,\ldots, g_{m-1})\;.$$
We can replace $\sigma_j$ with a $\sigma$, independent of $j$, defined by
$$\sigma=\max\{\sigma_j\ \vert\ fd\rho_j\neq0\}\;.$$
So the last sum is equal to
$$\sum T\llcorner (\sigma_j,\rho_j)(f,g_1,\ldots, g_{m-1})=T(\sigma f, \sum\rho_j, g_1,\ldots, g_{m-1})=0$$
 because $\sum\rho_j$ is constantly equal to $1$ on the support of $\sigma f$ (which is the support of $f$).

Therefore $d\circ R_{m-1}^j=R_{m}^j\circ d$ and then obviously $d\circ C_{m-1}=C_m\circ d$.

\smallskip

\noindent{\emph{(2) }} Upon taking a refinement of our open cover, we can suppose that there exists $c_1$, a positive constant, such that $ \Lip(\phi_j), \Lip(\phi_j^{-1})\leq c_1$. Let us denote by $\mu$ the (metric) mass measure of $T$ and by $m_j$ the (classical) mass measure of $R_m^j(T)$; we have
$$\|T\llcorner\rho_j\|(B)=(\mu\cdot\rho_j)(B)\;,$$
for every Borel set $B$, and
$$(\phi_j)_\sharp(\mu\cdot \rho_j)(B')\leq c_1(\mu\cdot\rho_j)(\phi_j^{-1}(B'))\;,$$
so
$$m_j(B)\leq c_1^2{N\choose m}(\mu\cdot\rho_j)(B)\;.$$
Summing on $j$ and denoting by $\mu'$ the mass measure of $C_m(T)$, we obtain
$$\mu'(B)\leq c_1^2{N\choose m}\mu(B)\;.$$
To obtain the other estimate, we note that
$$C_m(T)\llcorner \rho_j=R_m^j(T)$$
because
$$C_m(T)\llcorner\rho_j(fdg_1\wedge\ldots\wedge dg_m)=C_m(T)(f\rho_jdg_1\wedge\ldots\wedge dg_m)=$$
$$T(f\rho_j, g_1,\ldots, g_m)=R_m^j(T)(fdg_1\wedge\ldots\wedge dg_m)\;;$$
therefore
$$C_m^{-1}=\sum_j(\phi_j)_\sharp\circ (C_m^j)^{-1}\circ(\phi_j^{-1})_*(T\llcorner\rho_j)$$
with $T$ a classical current. So we obtain the estimate
$$\mu(B)\leq c_1^2\mu'(B)\;.$$

\smallskip

\noindent{\emph{(3)} and \emph{(4) }} The class of locally normal currents is stable under pushforward and contraction by a smooth function; the same is true for locally flat currents. Therefore these two conclusions follow easily from the corresponding ones in \ref{teo_comp}. \enddemo

\section{Bidimension}

Let $X$ be a (finite dimensional) complex space and $\Ol_X$ be the sheaf of holomorphic functions on $X$.

\begin{defin}\label{def_pq}Given $T\in\Di_k(X)$ (or $T\in M_k(X)$), we say that $T$ is of \emph{bidimension} $(p,q)$, with $p+q=k$, if
$$T(f,\pi_1,\ldots,\pi_k)=0$$
whenever $\exists\ J\subset\{1,\ldots, k\}$, with $|J|=p+1$, such that $\pi_j\in\Ol_X(\Omega)$ for every $j\in\ J$, for some open set $\Omega$ containing $\supp f$ or whenever $\exists\ J\subset\{1,\ldots, k\}$, with $|J|=q+1$, such that $\overline{\pi_j}\in\Ol_X(\Omega)$ for every $j\in J$, for some open set $\Omega$ containing $\supp f$.\end{defin}

We denote the space of local metric $(p,q)-$currents by $\Di_{p,q}(X)$; we also introduce the spaces $M_{p,q}(X)$, $M_{(p,q),\loc}(X)$, $N_{p,q}(X)$, $N_{(p,q),\loc}(X)$, their definition being obvious.

\begin{ex}If $X$ is a complex space with $\dim_\C X_\rg=n$, then the current of integration on the regular part of $X$ is an example of local metric $(n,n)$-current on $X_\rg$; we denote such a current by $[X]$ and we know, by \cite{Lelong}, that $[X]$ extends to an element of $N_{(n,n),\loc}(X)$, with $d[X]=0$.\end{ex}

The following Proposition summarizes some of the most important properties of $(p,q)-$currents.

\begin{propos}\label{prp_prop_pq}Let $p,\ q,\ r,\ s$ be non-negative integers such that $p+q=r+s=k$. Then
\begin{enumerate}\item $\Di_{p,q}(X)\cap\Di_{r,s}(X)=\emptyset$ if $(p,q)\neq(r,s)$;
\item if $\phi:X\to Y$ is holomorphic and proper and $T\in\Di_{p,q}(X)$, then $\phi_\sharp T\in \Di_{p,q}(Y)$;
\item $\Di_{p,q}(X)=\emptyset$ if $p>n$ or $q>n$ (with $n=\dim_\C X_\rg$);
\item if $T_j\to T\in\Di_k(X)$ weakly and $T_j\in\Di_{p,q}(X)$ for every $j$, then $T\in\Di_{p,q}(X)$;
\item $M_{p,q}(X)$ is closed in the Banach space $(M_k(X), \|\cdot\|)$.
\end{enumerate}
\end{propos}
\demo 
We note that the elements of the unitary $*-$algebra generated by $\Ol_X$ are dense in the germs of continuous functions, by the Stone-Weierstrass theorem; moreover, if $\pi$ belongs to that $*-$algebra, we can write it as $h_1\overline{g_1}+\ldots+h_m\overline{g_m}$, with $h_j,\ g_j\in\Ol_X$, therefore, by the product rule and multilinearity of currents, to evaluate $T(fd\pi)$, with $\pi_1,\ldots, \pi_k$ belonging to that algebra, it is enough to know $T(fd\pi)$ when each $\pi_j$ is either holomorphic or antiholomorphic.

\medskip

In particular, if $T(fd\pi)=0$ for every choice of $\pi_j$ in $\Ol_X$ or in $\overline{\Ol}_X$, then $T(fd\pi)=0$ for every $k-$tuple $\pi$; this implies \emph{(1)}. 

\medskip

To show \emph{(2)}, it is enough to notice that the composition with a holomorphic map preserves both holomorphic and antiholomorphic functions.

\medskip

Let us suppose that $p>n$ (the case $q>n$ is identical) and, initially, let us consider a metric form $fd\pi$ with each $\pi_j$ either holomorphic or antiholomorphic and $\supp f\subset V\Subset X_\rg$, where $V$ is some complex coordinate chart. Then, if $T\in \Di_{p,q}(X)$, $T(fd\pi)$ will vanish unless $p>n$ of the $\pi_j'$s are holomorphic, but then, by the chain rule, $T(fd\pi)=0$, because there are only $n$ holomorphic coordinate functions on $V$. This implies that $T$ is zero on every form compactly supported in $X_\rg$, then $\supp T\subset X_\sg$ and therefore $T=i_\sharp S$, with $S\in \Di_{p,q}(X_\sg)$ and $i:X_\sg\to X$ the inclusion (see Proposition \ref{prp_supp} below). Now \emph{(3)} follows by induction on the dimension.

\medskip

Property \emph{(4)} follows directly by the definition of $(p,q)-$current; moreover, convergence in mass implies weak convergence, thus \emph{(5)} follows as well. \enddemo

Given a (local) metric $k-$current $T$, we will say that $T$ \emph{admits a Dolbeault decomposition} if there exists $T_{p,q}$ of bidimension $(p,q)$, for each $(p,q)$ such that $p+q=k$, so that
$$T=\sum_{p+q=k}T_{p,q}\;.$$
The current $T_{p,q}$ will be called a \emph{$(p,q)-$component} of $T$.

\begin{propos}\label{prp_unique}The Dolbeault decomposition of $T$, if it exists, is unique.\end{propos}
\demo Let
$$T=\sum_{p+q=k}T_{p,q}=\sum_{p+q=k} S_{p,q}\;;$$
if $fd\pi$ is a metric form with $p$ holomorphic differentials and $q$ antiholomorphic differentials, then
$$T(fd\pi)=T_{p,q}(fd\pi)=S_{p,q}(fd\pi)\;.$$
Therefore, by the reasoning explained at the beginning of the proof of the previous Proposition, $T_{p,q}$ and $S_{p,q}$ coincide on the forms whose differentials belong to the unitary $*-$algebra generated by $\Ol_X$, so they coincide as currents.\enddemo

An interesting, but difficult problem is to know whether a (local) metric current admits a Dolbeault decomposition. In order to investigate some examples, we recall an important property of metric currents.

\begin{propos}[\cites{AK,Lang}]\label{prp_supp}
Suppose that $T\in\Di_k(X)$ (or $M_k(X)$), then
$$T(f,\pi_1,\ldots,\pi_k)=0$$
if $f\vert_{\supp T}=0$ or, if $k\geq1$, whenever some $\pi_j$ is constant on $\{f\neq0\}\cap\supp T$.
\end{propos}

We also remark that, in $\R^n$, classical (locally) flat currents\footnote{that is, a flat current in the sense of Federer and Fleming, see \cite{Federer}} can be extended to (local) metric currents.

\begin{ex}Let $T\in M_1(\C)$ be the metric current defined by
$$T(f,\pi)=\int_{S^1}fd\pi\;.$$
When $fd\pi$ is a smooth compactly supported $1-$form, this integral defines a classical flat current, therefore it can be extended to a metric current, which turns out to be of finite mass; $T$ can be written as
$$T(f,\pi)=\int_{S^1}\frac{i}{2}\pair{fd\pi, z\de_z-\bar{z}\de_{\bar{z}}}d\H^1$$
so
$$T(f,\pi)=\int_{S^1}\frac{i}{2}\pair{fd\pi, z\de_z}d\H^1-\int_{S^1}\frac{i}{2}\pair{fd\pi, \bar{z}\de_{\bar{z}}}d\H^1$$
$$=T_{1,0}(fd\pi)+T_{0,1}(fd\pi)\;.$$
Now, let us consider the two metric forms $(f,z)$ and $(f,z^2\bar{z})$, which coincide on $S^1=\supp T=\supp T_{1,0}=\supp T_{0,1}$; we have
$$T_{1,0}(f,z)=\frac{i}{2}\int_{S^1}zf(z)d\H^1$$
$$T_{1,0}(f,z^2\bar{z})=i\int_{S^1}zf(z)d\H^1\;.$$
Therefore, the conclusions of Proposition \ref{prp_supp} don't hold for $T_{1,0}$, implying that it cannot be extended to a metric current.
\end{ex}

\begin{ex}Let $\Leb^2$ be the usual Lebesgue measure on the complex plane and let $T\in M_{2,\loc}(\C)$ be the local metric current given by
$$T(f,\pi_1,\pi_2)=\int_{\C}\frac{1}{z}f\det(\nabla\pi)d\Leb^2\;.$$
We know that $T$ is a local metric current by \cite[Theorem 2.6]{Lang}; let $S=dT$ be its boundary. We compute the classical $(0,1)-$component of $S$, obtaining that 
$$S_{0,1}(f,\pi)=Cf(0)\frac{\de\pi}{\de z}(0)$$
for some constant $C$. It is easy to see that $S_{0,1}$ doesn't fullfill the conclusions of Proposition \ref{prp_supp}.
\end{ex}

\begin{ex}Let $Y$ be a $k-dimensional$ complex analytic set in $X$ and suppose $g\in L^1_{\loc}(Y,\H^{2k})$; then the current $T=[Y]\llcorner g (f,\pi)$ is of bidimension $(k,k)$. In general, if $T\in M_{k,\loc}(X)$, $\supp T\cap X_\sg$ is $\|T\|-$negligible and the Dolbeault components of the associated vector-field $\vec{\mathsf{T}}_{p,q}(x)$ are tangent to $\supp T$ for $\|T\|-$a.e. $x\in\supp T$, then $T$ admits a Dolbeault decomposition.\end{ex}

\subsection{The operators $\partial$ and $\debar$}

Let us suppose that $T$ is of bidimension $(p,q)$, with $p+q=k$, and that its boundary admits a Dolbeault decomposition. Then we can define $\de T$ and $\debar T$ as follows.

Let us write $dT=S_1+\ldots+S_h$ with $S_i$ of bidimension $(p_i,q_i)$; necessarily, $p_i+q_i=p+q-1=m$. If $(f,\pi)$ is a metric form with $p_i$ holomorphic differentials and $q_i$ antiholomorphic differentials, then
$$S_{p_i,q_i}(f,\pi)=dT(f,\pi)=T(1,f,\pi)\;;$$
if $p_i>p$ or $q_i>q$, then $T(1,f,\pi)=0$, as $T$ is a $(p,q)-$current. Therefore, we can only have two cases: $p=p_i$ and $q-1=q_i$ or $p-1=p_i$ and $q=q_i$; so, we have
$$dT=S_{p,q-1}+S_{p-1,q}$$
and we put
$$\de T=S_{p-1,q}\qquad \debar T=S_{p,q-1}\;.$$
Therefore, given a current $U$ which admits a decomposition in $(p,q)$ components, whose boundaries are decomposable too, we can define $\de U$ and $\debar U$ as
$$\de U=\sum_{i=1}^h\de U_i\qquad \debar U=\sum_{i=1}^h\debar U_i$$
where $U=U_1+\ldots+U_h$ is the $(p,q)$ decomposition.

\medskip

\begin{propos}\label{prp_pushf_debar}
If $\phi:X\to Y$ is a proper holomorphic map between complex spaces, then, for every current $T\in \Di_k(X)$ for which $\de T$ and $\debar T$ are defined, the following hold:
$$\phi_\sharp \de T=\de \phi_\sharp T\qquad \phi_\sharp\debar T=\debar \phi_\sharp T\;.$$
\end{propos}
\demo
By Proposition \ref{prp_prop_pq}, the pushforward of a $(p,q)-$current is a $(p,q)-$current. So, if $T\in \Di_{p,q}(X)$ and $dT=S_{p,q-1}+S_{p-1,q}$, then 
$$d\phi_\sharp T=\phi_\sharp dT=\phi_\sharp S_{p,q-1}+\phi_\sharp S_{p-1,q}\;.$$
By Proposition \ref{prp_unique}, this is the Dolbeault decomposition of $d\phi_\sharp T$; by definition
$$\de\phi_\sharp T=(d\phi_\sharp T)_{p-1,q}=\phi_\sharp S_{p-1,q}=\phi_\sharp \de T$$
$$\debar\phi_\sharp T=(d\phi_\sharp T)_{p,q-1}=\phi_\sharp S_{p,q-1}=\phi_\sharp \debar T\;.$$
\enddemo

\medskip

Moreover, it is easy to check that $C_m\circ\de=\de\circ C_m$ and $C_m\circ\debar=\debar\circ C_m$, where $C_m$ is the map given by Theorems \ref{teo_comp}, \ref{teo_comp_man}.

\begin{propos}\label{prp_formule_de}
We have that $\de^2=\debar^2=0$ and $\de\debar=-\debar\de$.
\end{propos}
\demo
By the locality property, we have that $d^2T=0$; therefore
$$0=(\de+\debar)(\de+\debar)T=\de^2 T+(\de\debar T+\debar\de T)+\debar^2 T$$
and, as the right hand side is a decomposition in $(p,q)$ components, every term has to be zero. Therefore
$$\de^2T=0\qquad \debar^2 T=0\qquad\de\debar T+\debar\de T=0$$
and, as we didn't make any assumption on $T$, the thesis follows.\enddemo

\medskip

We give a formula for $\debar T$, for analytic subsets of $\C^n$.

\begin{lemma}\label{lmm_eq_debar} Suppose $X$ is an analytic subset of some open set $\Omega\subseteq \C^n$; given $T\in \Di_m(X)$ such that $\debar T$ exists as a metric current and $(f,g)\in\Di^{m-1}(X)$ with $\Ci^1$ coefficients, we have 
\begin{equation}\label{eq_debar}\debar T(f,g)=\sum_{j=1}^nT\left(\frac{\de f}{\de\bar{z}_j}, \bar{z}_j, g\right)+\sum_{k=1}^{m-1}(-1)^{k}\sum_{j=1}^n T\left(f, \bar{z}_j,\ldots, \frac{\de g_k}{\de\bar{z}_j},\ldots\right)\;.\end{equation}
\end{lemma}
\demo The formula clearly holds for a classical current, by integration by parts. We noted before that $C_{m-1}\circ\debar=\debar\circ C_m$ and we know that there exists a metric current $S=\debar T$; therefore the thesis follows. \enddemo

\medskip

We can define a multilinear, local functional of finite mass by (\ref{eq_debar}) and denote it by $\debar T$. We define the space
$$W_m(X)=\{T\in \Di_m(X)\ :\ \debar T\in\Di_{m-1}(X)\}$$
and we note that $\debar W_m(X)\subseteq W_{m-1}(X)$, by Proposition \ref{prp_formule_de}. The space $W_{p,q}(X)$ is defined in the same way; unfortunately, we don't have any decomposition theorem for $W_m$ in terms of $W_{p,q}$.

\medskip

\begin{propos} \label{prp_int_parti}Suppose $U\subset X$ can be embedded as an analytic set into $\C^n$. Given $T\in W_m(U)$, $(u,v)\in\Di^k(U)$ with $\Ci^1(U)$ coefficients, we have that
$$\debar (T\llcorner(u,v))=(-1)^{k}(\debar T)\llcorner(u,v)-(-1)^k \sum_{j=1}^m T\llcorner(\de u/\de \bar{z}_j, \bar{z}_j, v)-$$
$$-(-1)^k\sum_{h=1}^k(-1)^h\sum_{j=1}^nT(u,\bar{z}_j,\ldots, \de v_h/\de\bar{z}_j,\ldots)$$
where $z_1,\ldots, z_n$ are the coordinates of the embedding in $\C^n$.
\end{propos}
\demo By Lemma \ref{lmm_eq_debar}, we have
$$((\debar T)\llcorner(u,v))(f,g)=(\debar T)(uf,v,g)=\sum_{j=1}^nT\left(\frac{\de uf}{\de \bar{z}_j}, \bar{z}_j, v,g\right)+$$
$$+\sum_{h=1}^k(-1)^h\sum T\left(uf, \bar{z}_j,\ldots, \frac{\de v_h}{\de\bar{z}_j},\ldots, g\right)+$$
$$+\sum_{h=1}^{m-k}(-1)^{h+k}\sum_{j=1}^nT\left(uf, \bar{z}_j, v, \ldots, \frac{\de g_{h}}{\de\bar{z}_j},\ldots\right)$$
for every $(f,g)\in\Di^{m-k}(U)$ with $\Ci^1$ coefficients.

We have
$$T\left(\frac{\de uf}{\de \bar{z}_j}, \bar{z}_j, v,g\right)=T\left(u\frac{\de f}{\de\bar{z}_j}, \bar{z}_j, v,g\right)+T\left(f\frac{\de u}{\de \bar{z}_j}, \bar{z}_j, v,g\right)$$
and we notice that
$$\sum_{j=1}^nT\left(u\frac{\de f}{\de\bar{z}_j}, \bar{z}_j, v,g\right)=(-1)^k\sum_{j=1}^nT\left(u\frac{\de f}{\de\bar{z}_j}, v,\bar{z}_j,g\right)=$$
$$=(-1)^k\sum(T\llcorner(u,v))\left(\frac{\de f}{\de\bar{z}_j},\bar{z}_j,g\right)\;.$$
So, again by Lemma \ref{lmm_eq_debar}, we obtain
$$\sum_{j=1}^nT\left(u\frac{\de f}{\de\bar{z}_j}, \bar{z}_j, v,g\right)+\sum_{h=1}^{m-k}(-1)^{h+k}\sum_{j=1}^nT\left(uf, \bar{z}_j, v, \ldots, \frac{\de g_{h}}{\de\bar{z}_j},\ldots\right)=$$
$$=(-1)^k\debar (T\llcorner(u,v))(f,g)\;.$$

Therefore we have
$$\debar (T\llcorner(u,v))(f,g)=$$
$$(-1)^{k}(\debar T)\llcorner(u,v)(f,g)-(-1)^k \sum_{j=1}^n T\llcorner(\de u/\de \bar{z}_j, \bar{z}_j, v)(f,g)-$$
$$(-1)^k\sum_{h=1}^k(-1)^h\sum_{j=1}^nT(u,\bar{z}_j,\ldots, \de v_h/\de\bar{z}_j,\ldots)(f,g)\;.$$
Now, the algebra of $\Ci^1$ functions is dense in the algebra of locally Lipschitz functions, therefore the current $(\debar T)\llcorner(u,v)$ is uniquely determined by this formula, which therefore holds for every $(f,g)$ in $\Di^{m-k}(U)$. \enddemo

\medskip

Obviously, we have a similar formula for $\de T$.

\section{Completely reducible spaces}

Let $X=L_1\cup\ldots\cup L_m$ be the union of linear subspaces $L_i$ of $\C^n$, with dimension $k_i$, such that $L_i\not\subseteq L_j$ whenever $i\neq j$. Obviously, $X$ is an analytic subset of $\C^n$.

Let us denote by $X_1$ the singular set of $X$ and  let us suppose that we have indexed the subspaces such that $\dim L_i\geq\dim L_{i+1}$ for every $i$; we will also suppose that $X$ isn't contained in any proper subspace of $\C^n$. Now, we consider a set $\mathcal{B}=\{L_1\}\cup\{L_{i_1},\ldots, L_{i_k}\}$ such that
$$\bigoplus_{L\in\mathcal{B}}L=\C^n\qquad L\not\subseteq\bigoplus_{\substack{L'\in\mathcal{B}\\L\neq L'}}L'\quad \forall\ L\in\mathcal{B}\;.$$
We have the projections 
$$\pi_1:X\to L_1\qquad\mathrm{and}\qquad \pi_{i_h}:X\to L_{i_h}$$
and the inclusions $j_1$ and $j_{i_h}$. We define
$$S=\bigcup_{L\in\mathcal{B}}L\qquad S'=\bigcup_{L\not\in\mathcal{B}}L\;.$$

Given $T\in M_m(X)$, we consider the currents
$$T_1=(j_1\circ \pi_1)_\sharp T$$
$$T_{i_h}=(j_{i_h}\circ\pi_{i_h})_\sharp T$$
and the difference
$$R=T-T_1-T_{i_1}-\ldots-T_{i_k}\;.$$

Let us choose $(f,\xi)\in\E^m(X)$ such that $\supp(f)\subseteq (L_1\cup L_{i_1}\cup\ldots\cup L_{i_k})\setminus X_1$; then, we can find $f_L\in \Lip_b(L)$ for $L\in\mathcal{B}$ such that 
$$f=\sum_{L\in\mathcal{B}}f_L\;.$$

Therefore we have
$$T(f,\xi)=\sum_{L\in\mathcal{B}} T(f_L,\xi)=\sum_{L\in\mathcal{B}} T\llcorner\chi_{L\setminus X_1}(f_L,\xi)$$ $$=T_1\llcorner\chi_{L_1\setminus X_1}(f_{L_1},\xi)+\sum_{h=1}^k T_{i_h}\llcorner\chi_{L_{i_h}\setminus X_1}(f_{L_{i_h}},\xi)=T_1(f_{L_1},\xi)+\sum_{h=1}^k T_{i_h}(f_{L_{i_h}},\xi)$$
$$=T_1(f,\xi)+\sum_{h=1}^k T_{i_h}(f,\xi)\;.$$
Then,
$$\supp(R)\subseteq \bigcup_{L\not\in\mathcal{B}}L\cup X_1=X_R\;.$$
Now, as $\dim L_1\geq \dim L_i$ for every $i$, we have that the maximum of the dimension of the irreducible components of $X_R$ is less or equal to $\dim L_1$; we can repeat our argument on $X_R$ and obtain a decomposition of $R$. Eventually, the remainder will have support contained in $X_1$, whose irreducible components have dimension strictly less than $\dim L_1$. By induction, we obtain the following

\begin{propos}\label{prp_str1}Given $T\in M_m(X)$, there exist proper holomorphic maps $h_i:\C^{k_i}\to X$, currents $Z_i\in M_m(\C^{k_i})$ for $i=1, \ldots, N$ such that
$$T=\sum_{i=1}^N T_i=\sum_{i=1}^N (h_i)_\sharp Z_i\;.$$
\end{propos}

\medskip

\begin{corol}\label{cor_deb1}Let us suppose that $T\in M_{p,q}(X)$ is a $\debar-$closed current; then there exists a $(p,q+1)-$current $U$ with locally finite mass such that $\debar U = T$ on $X$.\end{corol}
\demo
Let us write $T=\sum (h_j)_\sharp Z_j$, by the previous Proposition; as all the maps we used are holomorphic, we will also have that $Z_j$ is of bidimension $(p,q)$ and $\debar Z_j=0$ for every $j$. For every $j$, we can solve the equation $\debar V_j=Z_j$ on $\C^{k_j}$ with a local metric current $V_j$ (e.g. by convolution with the Cauchy kernel); therefore we can construct the local current $U=\sum(h_j)_\sharp V_j$ (as the maps $h_j$ are proper) and observe that
$$\debar U=\sum_{j=1}^N\debar (h_j)_\sharp V_j=\sum_{j=1}^n(h_j)_\sharp \debar V_j=\sum_{j=1}^N(h_j)_\sharp Z_j=T\;.$$
\enddemo

\bigskip

Now, if $X$ is a complex space which is locally biholomorphic to a union of linear subspaces of $\C^N$, e.g. if $X$ can locally be realized as a normal crossings divisor, we have the following results.

\begin{propos}\label{prp_loc_struct}Let $X$ be as said, then given $\Omega\subseteq X$ open, $T\in M_{m,\loc}(\Omega)$ and $x\in\Omega$, we can find an open set $\omega$ containing $x$, whose closure is contained in $\Omega$, holomorphic maps $h_i:V_i\to \Omega$,  $1\leq i\leq k$, with $V_i\subset\C^{n_i}$ and metric currents $T_i\in M_{m,\loc}(V_i)$ such that
$$T-\sum_{i=1}^k(h_i)_\sharp T_i$$
has support disjoint from $\overline{\omega}$.\end{propos}
\demo It is enough to choose an open set $U$ containing $x$, with $U\Subset\Omega$, such that there exists a biholomorphism of $U$ with a union of linear subspaces of $\C^N$ intersected with a ball, then we can apply the result of Proposition \ref{prp_str1} to the current $T\llcorner\sigma$, with $\sigma$ Lipschitz and supported in $U$. We now set $\omega=\mathrm{Int}(\supp\sigma)$. \enddemo

\medskip

\begin{teorema}\label{teo_loc_deb} Let $X$ be as said; given $\Omega\subseteq X$ open, $T\in M_{m,\loc}(\Omega)$ such that $\debar T=0$ and $x\in\Omega$, we can find an open set $\omega$ containing $x$, whose closure is contained in $\Omega$ and a current $S\in M_{m+1,\loc}(U)$ such that $T-\debar S$ has support disjoint from $\overline{\omega}$.\end{teorema}
\demo By Proposition \ref{prp_loc_struct} (and with the same notation), we find $\omega$ such that $T-\sum (h_i)_\sharp T_i$ has support disjoint from its closure; moreover, we know that $\debar T$ has support disjoint from $\omega$ if and only if $\debar (h_i)_\sharp T_i$ does. So, we can solve $\debar S_i=T_i$ in $V_i$, by convolution, and we set
$$S=\sum_{i=1}^k (h_i)_\sharp S_i\;.$$
By Proposition \ref{prp_pushf_debar}, we have that
$$\debar S(f,\pi)=\sum_{i=1}^k\debar((h_i)_\sharp S_i))(f,\pi)=\sum_{i=1}^k(h_i)_\sharp(\debar S_i)(f,\pi)=$$
$$\sum_{i=1}^k(h_i)_\sharp(T_i)(f,\pi)=T(f,\pi)$$
for every $(f,\pi)\in\Di^m(U)$ with $\supp f\subset\omega$. So, $T-\debar S$ has support disjoint from $\overline{\omega}$. \enddemo

\subsection{Holomorphic currents}

We now investigate the solutions of 
$$\debar T=0$$ 
when $T$ is a $(p,n)-$current on a complex space $X$, with $\dim_\C X_\rg=n$, which is locally biholomorphic to a union of linear subspaces of some $\C^N$; in the smooth case, these currents correspond to $(n-p,0)-$forms with holomorphic coefficients.

\begin{propos}\label{prp_loc_hol_curr}Let $X$ be as said and $\Omega\subseteq X$ be an open set, then, given a $\debar-$closed current $T\in M_{(p,n),\loc}(\Omega)$ and a point $x\in \Omega$, there exist an open set $\omega\subseteq\Omega$ containing $x$, a finite number of holomorphic maps $h_i:V_i\to\Omega$, with $V_i\subseteq\C^{n}$ open sets, and holomorphic $(n-p)-$forms $f_i\in\Omega^p(V_i)$ such that 
$$T-\sum_{i=1}^k(h_i)_\sharp[f_i]$$
has support disjoint from $\overline{\omega}$.\end{propos}
\demo By Proposition \ref{prp_loc_struct} (and with the same notation), we find $\omega$ such that $T-\sum (h_i)_\sharp T_i$ has support disjoint from its closure; moreover, we know that $\debar T$ has support disjoint from $\omega$ if and only if each $(h_i)_\sharp T_i$ does. As $\dim_\C X_\rg=n$, every $V_i$ is an open set in $\C^n$; we have that $\debar T_i$ is zero on $h_i^{-1}(\omega)$ so, as $T_i$ is of bidimension $(p,n)$, we can find a $(n-p,0)-$form $f_i$ with holomorphic coefficients such that $T_i=[f_i]$ in $h_i^{-1}(\omega)$.

By pushforward, we have the thesis. \enddemo

\medskip

We remark that the currents described in the previous Proposition are locally flat in any local affine embedding of $X$. Let $K_p$ be the kernel of $\debar: M_{(p,n),\loc}\to M_{(p,n-1),\loc}$.

\begin{propos}Let $X$ be as said, then $K_p=\pi_*\Omega^{n-p}_{X^\nu}$, where $\pi:X^\nu\to X$ is the normalization and $\Omega^{n-p}_{X^\nu}$ is the sheaf of holomorphic $(n-p)-$forms.\end{propos}
\demo Given $\Omega\subset X$ which is biholomorphic to a union of open neighborhoods of $0$ in $\C^n$, the preimage $\pi^{-1}(\Omega)\subset X^\nu$ is the disjoint union of these neighborhoods and $X^\nu$ is smooth. By Proposition \ref{prp_loc_hol_curr}, we know that $\debar-$closed $(p,n)-$currents on $\Omega$ are the holomorphic currents on $\pi^{-1}(\Omega)$. So the thesis follows. \enddemo

\medskip

We remark that $\Omega\subseteq X$ is Stein if and only if $\pi^{-1}(\Omega)$ is Stein, so
$$H^q(\Omega, K_p)=H^q(\pi^{-1}(\Omega), \Omega^{p}_{X^\nu})=0\;.$$

\medskip

We are now ready to give a global version of Proposition \ref{teo_loc_deb}.

\begin{teorema} Let $X$ be a Stein space with completely reducible singularities, $T\in M_{(n,n-1),\loc}(X)$ a $\debar-$closed metric current; then there exists a metric current $S$ of bidimension $(n,n)$ such that $\debar S=T$.\end{teorema}
\demo We consider an open cover $\{V_i\}_{i\in\N}$ of $X$ such that we can solve the Cauchy-Riemann equation on every $V_i$; we obtain a collection $\{S_i\}_{i\in\N}$ of metric currents with $S_i\in M_{(n,n),\loc}(\overline{V_i})$ and $\debar S_i=T$ in $M_{(n,n-1),\loc}(V_i)$.

On the sets $V_{ij}=V_i\cap V_j$, we have that $R_{ij}=S_i-S_j$ is a $\debar-$closed $(n,n)$ metric current in $M_{(n,n),\loc}(\overline{V_{ij}})$; lifting the cover $\{V_i\}$ to a cover $\{\Omega_j\}$ of the normalization $Y$ of $X$, we can find holomorphic functions $f_{ij}\in\Ol(\overline{\Omega_{ij}})$ such that $R_{\nu(i)\nu(j)}=\pi_\sharp[\Omega_{ij}]\llcorner f_{ij}$.

We now recall that the normalization of a Stein space is Stein, hence $H^1(Y,\Ol)=0$, therefore there exist functions $f_i\in\Ol(\Omega_i)$ such that $f_{ij}=f_i-f_j$. Defining $R_i=\pi_\sharp[\Omega_i]\llcorner f_i\in M_{(n,n),\loc}(\overline{V_i})$, on $\overline{V_{ij}}$ we have 
$$R_{ij}=R_i-R_j$$
so
$$S_i-R_i=S_j-R_j\;.$$
Thus we can define a metric current $S$ such that $\debar S=T$.\enddemo

\section{Complex Banach spaces}

In the following sections, we focus our attention on complex Banach spaces; let $E$ be such a space. The local version of metric currents developed by U. Lang doesn't exist in $E$, due to the fact that the compact sets in $E$ have empty interior; we will propose a replacement for these local objects in the next section.

Here we examine the behavior of metric currents in relation with their projections on finite dimensional subspaces of $E$; in order to recover informations on the whole space from its finite dimensional subspaces, we give the following definition (see also \cite{Noverraz}).

\begin{defin}\label{def_pap}A Banach space $E$ is said to have the \emph{projective approximation property (PA property)} if there exist a constant $a$ and an increasing collection $\{E_t\}_{t\in T}$ of finite dimensional subspaces of $E$ such that
\begin{enumerate}
\item $\{E_t\}_{t\in T}$ is a directed set for the inclusion;
\item $\displaystyle{E=\overline{\bigcup_{t\in T} E_t}}$;
\item for every $t\in T$ there exists a projection $p_t:E\to E_t$ with $\|p_t\|\leq a$.
\end{enumerate}
\end{defin}

Every Banach space with a Schauder basis has the PA property; more generally, the spaces $\Ci(K)$ of continuous functions on a compact set with the sup norm and $L^p(X,\mu)$, with $1\leq p\leq+\infty$ and $\mu$ a positive Radon measure on a locally compact space $X$. In this section, we will work with Banach spaces having the PA property; we will endow $T$ with the partial ordering coming from the inclusion relation between subspaces of $E$.

\medskip

\begin{propos}\label{prp_limlip}Let $f\in \Lip(E)$ and define $f_t=f\circ\pi_t$; then $f_t\to f$ pointwise and $\Lip(f_t)\leq a\Lip(f)$, for every $t$.
\end{propos}
\demo By property \emph{(2)} in Definition \ref{def_pap}, for every $x\in E$ there exists a sequence $\{x_j\}\subset E$, with $x_j\in E_{t_j}$, such that $x_j\to x$; by property \emph{(1)}, for a given $j$, we have that if $t>t_j$, then $f_t(x_j)=f(x_j)$. Moreover, 
$$\Lip(f_t)\leq \Lip(f)\cdot\|p_t\|\leq a\Lip(f)\;;$$
so, given $t, t'\in T$, let $j$ be such that $t_j\leq t$ and $t_j\leq t'$, then $f_t(x_j)=f_{t'}(x_j)=f(x_j)$ and
$$|f_t(x)-f_{t'}(x)|=|f_t(x)-f_t(x_j)+f_{t'}(x_j)-f{t'}(x_j)|\leq 2a\Lip(f)\|x-x_j\|$$
which goes to $0$ as $j\to\infty$. \enddemo

\begin{propos}Let $T\in M_k(E)$ and define $T_t=(\pi_t)_\sharp(T)\in M_k(E_t)$ for every $t\in T$ such that $\dim_{\C}E_t\geq k$; through the inclusion $i_t:E_t\to E$, we can look at $T_t$ as an element of $M_k(E)$ again. Then, $T_t\to T$ weakly.\end{propos}
\demo Let $\mu_t$ be the mass of $T_t$ and $\mu$ the mass of $T$; then $\mu_t=(\pi_t)_\sharp \mu$. By \cite[Lemma 2.9]{AK}, the support of $\mu$ is a $\sigma-$compact set, therefore for every $\epsilon>0$ there exists a compact $K_\epsilon$ such that $\mu(E\setminus K_\epsilon)\leq \epsilon$. As $\pi_t\to \mathrm{Id}_E$ uniformly on every compact set (because of the PA property), we have that $\mu_t\to \mu$ on $K_\epsilon$, which implies that
$$\int_E f\circ \pi_t d\mu=\int_{E}fd\mu_t\to \int_E fd\mu$$
for every $f\in \Lip_b(E)$. By this result and by Proposition \ref{prp_limlip}
$$T_t(f,\pi)=T(f\circ\pi_t,\pi\circ\pi_t)\to T(f,\pi)$$
which is our thesis. \enddemo

\medskip

\begin{defin}Let $\{E_t, p_t\}_{t\in T}$ be the countable collection of subspaces and projections given by the PA property. We call it a \emph{projective approximating sequence} (PAS) if $p_t\circ p_s=p_{\min\{s,t\}}$.\end{defin}

We note that every separable Hilbert space or, more generally, every Banach space with a Schauder basis contains a PAS.

\begin{teorema}\label{teo_ext}Let us suppose that $\{E_t, p_t\}$ is a PAS in $E$. If we are given a collection of metric currents $\{T_t\}_{t\in T}$ such that 
\begin{enumerate}
\item $T_t\in N_k(E_t)$,
\item $(p_{t}\vert_{E_{t'}})_\sharp T_{t'}= T_{t}$ for every $t, t'\in T$ with $t'>t$,
\item $\|T_t\|\leq (p_t)_*\mu$ and $\|d T_t\|\leq (p_t)_*\nu$ for every $t\in T$ and some $\mu,\ \nu$ finite Radon measures on $E$, 
\end{enumerate}
then there exists $T\in N_k(E)$ such that $(p_t)_\sharp T=T_t$ for every $t\in T$. \end{teorema}

\demo We consider, in $\E^k(E)$, the subspaces $(p_t)^*\E^k(E_t)$: their union $\mathcal{P}^k$ is dense, with respect to pointwise convergence, with bounded Lipschitz constants. We define a functional $T:\mathcal{P}^k\to\C$, by setting $T(f,\pi)=T_t(f,\pi)$, with $t$ such that $(f,\pi)\in(p_t)^*\E^k(E_t)$. By hypothesis \emph{(2)}, this definition is well posed; the functional so defined is obviously multilinear and local on $\mathcal{P}^k$, moreover, by hypothesis \emph{(3)}, we have that there exists a finite Radon measure $\mu$ on $E$ such that
$$T(f,\pi)\leq \prod_{j=1}^k \Lip(\pi_j) \int_E|f|d\mu\qquad \forall\ (f,\pi)\in\mathcal{P}^k\;.$$
In the same way, we know that there exists a finite Radon measure $\nu$ such that
$$dT(f,\pi)\leq\prod_{j=1}^{k-1} \Lip(\pi_j)\int_E|f|d\nu\qquad\forall\ (f,\pi)\in\mathcal{P}^{k-1}\;;$$
therefore, the current $T$ is also continuous on $\mathcal{P}^k$, being normal on this set of metric forms (see \cite[Proposition~5.1]{AK}). 

Extending $T$ by density, we obtain a multilinear, local, continuous functional on $\E^k(E)$, whose mass is bounded by $\mu$ and whose boundary's mass is bounded by $\nu$; thus, the extension is a normal current, which we denote again with $T$, and it is not hard to check that it verifies $(p_t)_\sharp T=T_t$ $\forall\ t\in T$. \enddemo

\medskip

We can substitute the request of the existence of a PAS and the compatibility condition (hypothesis \emph{(2)}) with an assumption on the existence of a global object. A \emph{metric functional} is a function $T:\E^k(E)\to\C$ which is subadditive and positively $1-$homogeneous with respect to every variable. For metric functionals, we can define mass, boundary and pushforward (see Section 2 of \cite{AK}).

\begin{propos}Let $E$ be a Banach space with the PA property; suppose that $T:\E^k(E)\to\C$ is a metric functional with finite mass, whose boundary has finite mass too, such that $(p_t)_\sharp T\in N_k(E_t)$ for every $t\in T$. Then $T\in N_k(E)$.\end{propos}
\demo We have 
$$\|(p_t)_\sharp T\|\leq a(p_t)_*\|T\|$$
$$\|d(p_t)_\sharp T\|\leq a(p_t)_*\|dT\|\;,$$ 
so, by  the previous Theorem there exists $\widetilde{T}\in N_k(E)$ such that $(p_t)_\sharp \widetilde{T}=T_t$. This means that $\widetilde{T}$ and $T$ coincide on the metric forms in $\mathcal{P}^k$; by density, we conclude that $\widetilde{T}=T$. \enddemo

\subsection{Bidimension}

Definition \ref{def_pq} is meaningful also when $X$ is a complex Banach space; for a careful analysis of the infinite dimensional holomorphy, we refer the interested reader to the first chapters in \cite{Noverraz}. Here we only notice that Lipschitz analytic functions are not necessarily dense in Lipschitz functions and we cannot work with local concepts as in Section 3, because the spaces of local currents do not make sense on a Banach space.

However, inspired by the links we have found between the finite dimensional projections of a current and the current itself, we would like to give a different characterization of $(p,q)-$currents. We say that $T\in M_{k}(E)$ is \emph{finitely} of bidimension $(p,q)$ if every finite dimensional projection of it is a $(p,q)-$current.

\begin{propos} $T\in M_k(E)$ is a $(p,q)-$current if and only if it is finitely so.\end{propos}
\demo A projection $p:E\to V$ is a continuous complex linear operator, thus holomorphic, so $T\in M_{p,q}(E)$ implies $p_\sharp T\in M_{p,q}(V)$, so one implication is proved.

On the other hand, if $h\in \Ol(E)$, then $h\vert_{E_t}\in \Ol(E_t)$ for every $t\in T$; so, if $(f,\pi)\in \E^k(E)$ contains $p+1$ holomorphic differentials, then so does $(f\vert_{E_t}, \pi\vert_{E_t})\in \E^k(E_t)$. Therefore,
$$T(f,\pi)=\lim_{t\in T} T(f\circ p_t, \pi\circ p_t)=\lim_{t\in T} T_{t}(f,\pi)\;.$$
As $T_t$ is a finite dimensional projection, it is of bidimension $(p,q)$, so the right hand side is zero. The same argument applies when $(f,\pi)$ contains $q+1$ antiholomorphic differentials, giving us the desired conclusion. \enddemo

As an application of Theorem \ref{teo_ext}, we have the following result about the existence of a Dolbeault decomposition for $T\in M_k(E)$.

\begin{propos}\label{prp_dec}Let us suppose that $\{E_t, p_t\}$ is a PAS in $E$. Let $T\in N_k(E)$; if $T_t$ has a Dolbeault decomposition in normal $(p,q)-$currents in $E_t$ for all $t\in T$, with a finite Radon measure $\nu$ (independent of $t$) whose pushforward dominates the boundaries' masses, then also $T$ admits a Dolbeault decomposition.\end{propos}
\demo Let us fix a pair $(p,q)$ such that $p+q=k$ and let $S_t$ be the $(p,q)-$component of $T_t$; by hypothesis, $S_t\in N_{p,q}(E_t)$ and $\|d S_t\|\leq (p_t)_*\nu$, independently of $t$, and it is not hard to show that $\|S_t\|\leq C'\|T_t\|\leq C''\|T\|$, with $C',\ C''$ independent of $t$ (in particular, independent of $\dim E_t$).

Last thing to check is the compatibility condition (condition \emph{(2)} in Theorem \ref{teo_ext}), but this follows easily from the invariance of the bidimension under pushforward by holomorphic maps. Applying Theorem \ref{teo_ext}, we have the thesis. \enddemo

\begin{rem} In general it isn't easy to verify the hypotheses of Proposition \ref{prp_dec} for a current $T\in N_k(E)$; however, this result is an example of a general phenomenon: in a Banach space with the projective approximation property, it is often enough to verify a certain property for finite dimensional subspaces in order to obtain that it holds for the whole space. For instance, any equality between currents holds in $E$ if and only if it holds finitely, that is, whenever the currents are pushed forward through a finite rank projection.\end{rem}

Employing the idea given in this Remark, we can show the following.

\begin{corol}If $T\in M_k(E)$ admits a Dolbeault decomposition, then it is unique.\end{corol}

\section{Quasi-local metric currents}

To partially overcome the problems related to the lack of a local theory of currents, we introduce a new definition, which somehow locates midways between the local and the global one.

Let $\E_{\q}(E)$ be defined as follows:
$$\E_{\q}(E)=\bigcup_{R>0}\{f\in\Lip(E)\ :\ \supp f\subset B(0,R)\}$$
where $B(0,R)$ is the ball of centre $0$ and radius $R$ in $E$. We say that a sequence $\{f_j\}\subset \E_{\q}(E)$ converges to $f\in \E_{\q}(E)$ if $f_j\to f$ pointwise, $\Lip(f_j), \Lip(f)\leq C$ and $\supp f_j, \supp f\subset B(0,R)$ for some $R$.

We also define $\Lip_{\q}(E)$ to be
$$\{ f\in \Ci^0(E)\ :\ f\vert_{B(0,R)}\in \Lip(B(0,R))\ \forall\ R>0\}\;.$$
A sequence $\{\pi_j\}$ in this space converges if it converges pointwise and the Lipschitz constants on any fixed ball are uniformly bounded in $j$ (but not necessarily with respect to the radius of the ball).

Finally, we define the spaces of quasi-local metric forms as
$$\E^k_\q(E)=\E_\q(E)\times[\Lip_\q(E)]^k\;.$$

\begin{defin}A \emph{quasi-local} $k-$dimensional metric current is a functional $T:\E^k_\q(E)\to\C$ satisfying the following
\begin{enumerate}
\item $T$ is multilinear
\item whenever $(f^i,\pi^i)\to (f,\pi)$ in $\E^k_\q(E)$, $T(f^i,\pi^i)\to T(f,\pi)$
\item $T(f,\pi)=0$ whenever there is an index $j$ such that $\pi_j$ is constant on a neighborhood of $\supp f$
\item for every $R>0$ there is a finite Radon measure $\mu_R$ such that
$$|T(f,\pi)|\leq\prod_{j=1}^k\Lip(\pi_j\vert_{B(0,r)})\int\limits_{B(0,R)}|f|d\mu_R$$
for every $(f,\pi)\in\E^k_\q(E)$ with $\supp f\subset B(0,R)$.\end{enumerate}
We denote the space of such currents by $M_{k,\q}(E)$.
\end{defin}

The last condition can be rephrased as: there exits a Radon measure $\mu$ on $E$, such that $\mu(B(0,R))<+\infty$ for every $R>0$ and such that
$$|T(f,\pi)|\leq\prod_{j=1}^k\Lip(\pi_j\vert_{\supp f})\int\limits_{\supp f}|f|d\mu$$
for every $(f,\pi)\in\E^k_\q(E)$. If $\mu$ happens to be a \emph{finite} Radon measure, then $T$ is indeed a $k-$dimensional metric current in the sense of Ambrosio and Kirchheim.

\medskip

The definitions of boundary, pushforward, contraction are the same of the usual metric currents; the pushforward can be performed with any \emph{quasi-local} proper Lipschitz map, that is any map which is Lipschitz on $B(0,R)$ for every $R>0$ and such that the preimage of any bounded set is bounded.

\begin{rem} The projections $p_t$ are by no means quasi-local, so we cannot repeat verbatim the arguments of the previous section.\end{rem}

The space $N_{k,\q}(E)$ is defined as the set of quasi-local currents whose boundary is again a quasi-local current, that is, has quasi-locally finite mass.

\medskip

By the mass condition, we can extend any $T\in M_{k,\q}(E)$ to a functional on $k+1-$tuples $(f,\pi)$ where $\pi\in [\Lip_{\textrm{q-loc}}(E)]^k$ and $f\in\mathcal{B}_b^{\infty}(E)$, that is the algebra of bounded Borel functions with bounded support in $E$. The basic properties of metric currents hold also for this quasi-local variant. Namely, we have the following.

\begin{propos}Given $T\in M_{k,\q}(E)$, we denote again by $T$ its extension to $\mathcal{B}^\infty_b(E)\times[\Lip_{\textrm{q-loc}}(E)]^k$; then
\begin{enumerate}
\item $T$ is multilinear in $(f,\pi)$ and
$$T(fd\pi_1\wedge\ldots\wedge d\pi_k)+T(\pi_1df\wedge\ldots\wedge d\pi_k)=T(\sigma d(f\pi_1)\wedge\ldots\wedge d\pi_k)$$
whenever $f,\pi\in \E_\q(E)$ and $\sigma\in\mathcal{B}^\infty_b(E)$ is equal to $1$ on the support of $f\pi_1$ and
$$T(fd\psi_1(\pi)\wedge\ldots\wedge d\psi_k(\pi))=T(f\det\nabla\psi(\pi)d\pi_1\wedge\ldots\wedge d\pi_k)$$
whenever $\psi=(\psi_1,\ldots,\psi_k)\in\mathcal{C}^1(\R^k,\R^k)$;
\item $$\lim_{i\to\infty}T(f^i,\pi^i_1,\ldots,\pi^i_k)=T(f,\pi)$$
whenever $f^i-f\to0$ in $L^1(E,\|T\|)$ and $\pi^i_j\to\pi_j$;
\item $T(f,\pi)=0$ if $\{f\neq 0\}=\bigcup B_i$ with $B_i\in\mathcal{B}(E)$ and $\pi_i$ constant on $B_i$.
\end{enumerate}
\end{propos}

The definition of $(p,q)-$current given in Section 3 can be applied also to quasi-local currents; we have the same results of Propositions \ref{prp_prop_pq}, \ref{prp_unique}. Given $(f,\pi)\in \E^k_\q(E)$, the differentials $\pi_1,\ldots, \pi_k$ can be approximated by analytic functions of a finite number of variables, so the proofs go on almost identically.

\section{Quasi-local solution to $\debar U=T$}

Given a function $a:E\to\C$, we set
$$a_t(x)=a(tx)\qquad \forall\ t\in\C\;.$$

Let $T\in N_{k}(E)$ be a $(0,k)-$current, with $\supp T$ bounded and $0\not\in\supp T$; we define the following $(k+1)-$dimensional metric functional
$$C_{\debar}(T)(f,\pi_1,\ldots,\pi_{k+1})=\frac{1}{2\pi i}\sum_{j=1}^{k+1}(-1)^{j+1}\int_{\C}T\left(f_t\frac{\de \pi_{jt}}{\de \bar{t}}d\widehat{\pi}_{jt}\right)\frac{dt\wedge d\bar{t}}{t-1}\;,$$
where $d\widehat{\pi}_j$ is the wedge product of all the differentials different from $\pi_j$.

\begin{lemma}$C_{\debar}(T)$ is a multilinear, local metric functional, with quasi-locally finite mass.\end{lemma}
\demo Multilinearity and locality are obvious. Let $B(0,r)$ be a ball containing the support of $T$ and let $B(0,d)$ be a ball disjoint from the support of $T$.

We have that, for $fd\pi\in\E^{k+1}_\q(E)$, with $\Lip(\pi_i)=1$ and $\supp f\subset B(0,R)$, the following holds
$$|C_{\debar}(T)(fd\pi)|\leq \frac{r(k+1)}{2\pi}\int_{|t|<R/d}\frac{|t|^k}{|t-1|}\int_{E}|f_t|d\|T\|dt\wedge d\bar{t}\;.$$
Moreover, given a bounded borel set $A$, 
$$\|C_{\debar}(T)\|(A)\leq \frac{r(k+1)}{2\pi}\int_{|t|<R/d}\frac{|t|^k}{|t-1|}\|T\|(A/t)dt\wedge d\bar{t}\;.$$
We want to estimate it for $A=B(0,R)$, with $R>r$. We split the integral in $t$ in two parts: one small ball around the origin and the rest of the ball of radius $R/d$. In the small ball around the origin, we have
$$\int_{|t|<\epsilon}\frac{|t|^k}{|t-1|}\|T\|(B(0,R/t))dt\wedge d\bar{t}\leq c_1\epsilon^{2+k}M(T)\;.$$
On the rest of the outer ball, we have
$$\int_{\epsilon<|t|<R/d}\frac{|t|^k}{|t-1}\|T\|(B(0,R/t))dt\wedge d\bar{t}\leq (R/d)^{k+2}M(T)\;.$$
So, letting $\epsilon\to0$, we get
$$M_{B(0,R)}(C_{\debar}(T))\leq \frac{r(k+1)}{2\pi}(R/d)^{k+2}M(T)\;.$$
Therefore, the mass of $C_{\debar}(T)$ is quasi-locally finite. \enddemo

\begin{propos}\label{prp_d_cone}We have
$$d C_{\debar}(T)=C_{\debar}(d T) + T$$
as quasi-local metric functionals.
\end{propos}
\demo Let $fd\pi$ be a quasi-local metric $k-$form such that $f$, $\pi_1,\ldots, \pi_{k}$ have Lipschitz derivatives; we define the function
$$\phi(t)=T\left( f_td\pi_t\right)$$
and we note that, as $T$ is a quasi-local current, thus continuous, 
$$\frac{\de\phi}{\de\bar{t}}=T\left(\frac{\de f_t}{\de\bar{t}}d\pi_t\right)+\sum_{j=1}^k(-1)^{j+1}T\left(fd\frac{\de\pi_{jt}}{\de\bar{t}}\wedge d\widehat{\pi}_{jt}\right)\;.$$
By the definition of boundary this expression can be rewritten as
\begin{equation}\label{eq_bdry}T\left(\frac{\de f_t}{\de\bar{t}}d\pi_t\right)+\sum_{j=1}^k(-1)^{j}\left[T\left(\frac{\de\pi_{jt}}{\de\bar{t}}df_t\wedge d\widehat{\pi}_{jt}\right)-dT\left(f_t\frac{\de\pi_{jt}}{\de\bar{t}}d\widehat{\pi}_{jt}\right)\right]\;,\end{equation}
exactly as in the proof of Proposition 10.2 in \cite{AK}.

Given a generic form $fd\pi\in\E^{k}_\q(E)$, the conclusion still holds: it is enough to approximate $f$ and $\pi_j$ by
$$f^\epsilon(x)=\int_{\C}f(sx)\rho_\epsilon(s)ds\wedge d\bar{s}\;,\qquad \pi_j^\epsilon(x)=\int_{\C}\pi_j(sx)\rho_\epsilon(s)ds\wedge d\bar{s}\;,$$
where $\rho_\epsilon$ are convolution kernels, compactly supported, $w^*-$converging to $\delta_1$.

By Fubini's theorem we have
$$\lim_{\epsilon\to0}\frac{\de f^{\epsilon}_t}{\de\bar{t}}(x)=\frac{\de f_t}{\de\bar{t}}(x)\;,\quad\lim_{\epsilon\to0}\frac{\de \pi_{jt}^\epsilon}{\de\bar{t}}(x)=\frac{\de\pi_{jt}}{\de\bar{t}}(x)$$
for $\|T\|+\|dT\|-$a.e. $x$, for $\mathcal{L}^2-$a.e. $t$; therefore the derivative with respect to $\bar{t}$ of $t\mapsto T(f^\epsilon_t d\pi^\epsilon_t)$ converges for a.e. $t\in\C$.

We notice that the supports of convolutions, for $\epsilon$ small enough are not significantly distant from the supports of the original functions. 

As $d(C_{\debar}(T))(fd\pi)+C_{\debar}(dT)(fd\pi)$ is equal to the integral of the expression in (\ref{eq_bdry}), multiplied by $(t-1)^{-1}$, over $\C$, we have that
$$d(C_{\debar}(T))(fd\pi)+C_{\debar}(dT)(fd\pi)=\int_{\C}\frac{\de\phi(t)}{\de\bar{t}}\frac{dt\wedge d\bar{t}}{t-1}=\phi(1)=T(fd\pi)\;,$$
because, as $\supp f$ is bounded and $\supp T$ has positive distance from $0$, the function $\phi(t)$ is compactly supported. \enddemo

\begin{corol}$C_{\debar}(T)$ is in $N_{k+1,\q}(E)$ and of bidimension $(0,k+1)$.\end{corol}
\demo Employing the previous Proposition, we repeat the proof of Proposition 10.2 in \cite{AK}, obtaining the continuity of $C_{\debar}(T)$ as a quasi-local metric current. An easy calculation shows that, given $fd\pi\in E^{k}_\q(E)$ and $h\in\Ol(E)$, we have
$$C_{\debar}(T)(fdh\wedge d\pi)=0$$
so $C_{\debar}(T)$ is of bidimension $(0,k+1)$.\enddemo

\begin{corol}We have $\debar C_{\debar}(T)=C_{\debar}(\debar T)+T$.\end{corol}
\demo For a $(0,q)-$current $S$, $dS=\debar S$.\enddemo

Summing up the previous results, we have the following quasi-local solution of the Cauchy-Riemann equation.

\begin{teorema}\label{teo_0q} Let $T\in N_{k}(E)$ be a $(0,k)-$current, with $\debar T=0$, such that $\supp T$ is bounded and at a positive distance from $0$; then there exists a quasi-local metric $(0,k+1)-$current $U$ such that $\debar U=T$.\end{teorema}
\demo If $dT=\debar T=0$, then, letting $U=C_{\debar}(T)$, we have
$$\debar U=C_{\debar}(\debar T)+T=T$$
as quasi-local currents. \enddemo

\begin{rem}We can control the mass of the solution $U$ on a ball by the mass of $T$ on that ball, with a constant depending only on the support of $T$, the dimension of $T$ and the radius of the ball.\end{rem}

This result can at once be extended to $(p,q)-$currents, if we add the request that $\debar T=\de T=0$.

\begin{corol} Let $T\in N_{k}(E)$ be a $(p,q)-$current, with $\debar T=\de T=0$, such that $\supp T$ is bounded and at a positive distance from $0$; then there exists a quasi-local metric $(p,q+1)-$current $U$ such that $\debar U=T$.\end{corol}
\demo We perform the same cone construction and we note that, by Proposition \ref{prp_d_cone}, $dC_{\debar}(T)=T$. Now, if $fd\pi\in E^{k}_\q(E)$ has $p+1$ holomorphic differentials, say $\pi_1,\ldots, \pi_{p+1}$, and $h\in \Ol(E)$, then 
$$T\left(f_t\frac{\de h_t}{\de\bar{t}}d\pi\right)=T(0d\pi)=0$$
$$T\left(f_t\frac{\de \pi_{jt}}{\de\bar{t}}d\widehat{\pi}_{jt}\right)=T(0d\widehat{\pi}_{jt})=0\qquad\textrm{ for } 1\leq j\leq p+1$$
and $d\widehat{\pi}_{jt}$, for $j\geq p+2$, contains $p+1$ holomorphic differentials, so $T(\sigma d\widehat{\pi}_{jt})=0$ for any $\sigma$ with bounded support. This means that $C_{\debar}(T)$ is a $(p,q+1)-$current; as $dC_{\debar}(T)=T$ is of bidimension $(p,q)$, we conclude that $dC_{\debar}(T)=\debar C_{\debar}(T)$.\enddemo

However, with some more effort, we can obtain the general result for $(p,q)-$currents.

\begin{teorema}Let $T\in N_k(E)$ be a $(p,q)-$current, with $\debar T=0$, such that $\supp T$ is bounded and at a positive distance from $0$; then there exists a quasi local metric $(p,q+1)-$current $U$ such that $\debar U=T$. \end{teorema}
\demo We remark that, as $T$ is normal and $\debar T=0$, then $dT$ admits a Dolbeault decomposition, where $(dT)_{p-1,q}=\de T= dT$. Let $h_1,\ldots, h_p\in\Ol(E)$ be holomorphic functions and set $H=(h_1,\ldots, h_p)$; then $$S_H=T\llcorner(1,h_1,\ldots, h_p)$$
is a $(0,q)-$current such that 
$$\debar S_H=dS_H=(dT)\llcorner(1,h_1,\ldots, h_p)=(dT)_{p-1,q}\llcorner(1,h_1,\ldots, h_p)$$
and the last term is $0$ by the definition of $(p-1,q)-$current. Therefore $\debar S_H=0$.

Now, by Theorem \ref{teo_0q}, there exists a $(0,q+1)-$current $V_H$ such that $\debar V_H=S_H$;  for each $(f,H,\pi)\in\E^{k+1}_{\q}(E)$ with $H$ a $p-$tuple of holomorphic functions, we define the metric functional
$$U(f,H,\pi)=V_H(f,\pi)$$
and we set
$$U(f,\eta)=0$$
whenever $\eta$ contains at most $p-1$ holomorphic functions and at least $q+2$ antiholomorhpic functions.

It is easy to check that $U$ is then defined on every (quasi-local) $k+1-$form with either holomorphic or antiholomorphic differentials; this allows us to extend $U$ as a multilinear, local, alternating functional on the (quasi-local) $k+1-$forms with analytic coefficients.

We have, whenever $\supp f\subset B(0,R)$,
$$|U(f,H,\pi)|=|V_H(f,\pi)|\leq\prod_{j=1}^{q+1}\Lip(\pi_j)C(R)\int_{B(0,R)}|f|d\|S_H\|\leq$$
$$\leq \prod_{j=1}^p\Lip(h_j)\prod_{j=1}^q\Lip(\pi_j)C(R)\int_{B(0,R)}|f|d\|T\|\;.$$
So $\|U\|_{B(0,R)}\leq C(R)\|T\|_{B(0,R)}$, which means that the mass of $U$ is quasi-locally finite, wherever $U$ is defined.

Moreover, let $(f,\eta)$ be a $k-$form with each of $\eta_1,\ldots,\eta_k$ holomorphic or antiholomorphic and
$$f=\chi_{B(0,R)}\cdot g$$
for some $r<R$ and $g$ holomorphic or antiholomorphic. Then
$$dU(f,\eta)=U(\chi_{B(0,R)}, g,\eta)\;.$$
If $(g,\eta)$ contains more that $p$ holomorphic functions or more than $q+1$ antiholomorphic functions, then $dU(f,\eta)=0$.

If $\eta=(H,\pi)$, with $H=(h_1,\ldots, h_p)$ holomorphic, then
$$|dU(f,\eta)|=|U(\chi_{B(0,R)}, g,H,\pi)|=|V_H(\chi_{B(0,R)}, g,\pi)|=|(dV_H)(f,\pi)|$$
$$=|S_H(f,\pi)|=|T(f,H,\pi)|=|T(f,\eta)|$$
so, in this case $\|dU\|\leq\|T\|$.

If $\eta$ contains only $p-1$ holomorphic functions, then $g$ has to be holomorphic (otherwise $U(\chi_{B(0,R)},g,\eta)=0$); so we set $\eta=(\pi',H')$ and
$$|U(\chi_{B(0,R)},g,\pi',H')|=|V_{gH'}(\chi_{B(0,R)}, \pi')|$$
$$=|(dV_{gH'})(\pi'_1\chi_{B(0,R)}, \pi'_2,\ldots,\pi'_{q+1})|$$
$$=|S_{gH'}(\pi_1'\chi_{B(0,R)}, \pi_2',\ldots,\pi_{q+1}')|=|T(\pi'_1\chi_{B(0,R)}, g, H', \widehat{\pi'}_1)|$$
$$=|dT(f\pi_1', H', \widehat{\pi'}_1)-T(f, \pi_1',H', \widehat{\pi'}_1)|=|dT(f\pi_1', H', \widehat{\pi'}_1)|$$
where we have employed the definition of boundary and the definition of $(p,q)-$current. We can suppose, without loss of generality, that $\pi'_1(0)=0$, so that we get
$$\|dU\|\leq R\|dT\|\;.$$

\medskip

Summing up, we have that $\|dU\|\leq \max\{\|T\|, R\|dT\|\}$, on forms $(f,\eta)$ where every component is either holomorphic or antiholomorphic on the ball $B(0,R)\supset\supp f$.

This allows us to extend $U$ as a quasi-local current on forms $(f,\eta)\in\E^{k+1}_{\q}(E)$ with $f=g\chi_{B(0,R)}$ for some $R$ and $g$ holomorphic or antiholomorphic; by linearity in the first component, we can allow $g$ to be analytic, then by density (and quasi-local finiteness of mass) we can extend $U$ to $\E^k_{\q}(E)$.

By the previous computations, $U$ is then a quasi-locally normal current, of bidimension $(p,q+1)$, with $\debar U=T$; moreover, the mass of $U$ is controlled, on a fixed ball, by the mass of $T$ and we have $\|U\|_{B(0,R)}+\|dU\|_{B(0,R)}\leq A(R)\|T\|+B(R)\|dT\|$.
\enddemo

\begin{rem} The hypothesis that $0$ has positive distance from the support of $T$ can be avoided, by constructing the cone from a point different from the origin, as long as the support of $T$ is bounded.\end{rem}

The hypothesis of the boundedness of $\supp T$ seems much harder to get rid of and, to date, we do not even know if it is possible. In the same way, it is not apparent that one can improve the estimates on mass in order to obtain a metric current (not a quasi-local one) from the cone construction; in the one-dimensional case, this cone construction consists in the convolution with the Cauchy kernel, which does not in general give a compactly supported solution.

\medskip

Two natural questions arise:
\begin{enumerate}
\item are there conditions on $T$ which ensure that the solution obtained with the cone construction to $\debar U=T$  has bounded support? or finite mass?
\item can the alleged solution with bounded support be obtained as a minimizer for the mass or the quasi-local mass among all quasi-local solutions to the Cauchy-Riemann equation?
\end{enumerate}

\bibliographystyle{amsalpha}
\bibliography{some_app}\end{document}